\documentclass[11pt]{article}
\usepackage{color,amsfonts,amssymb}
\usepackage{amsfonts,epsf,amsmath}
\usepackage{latexsym}
\usepackage{graphicx}
\usepackage{array,float}

\newtheorem{theorem}{Theorem}[section]
\newtheorem{lemma}[theorem]{Lemma}
\newtheorem{proposition}[theorem]{Proposition}

\newcommand{\proof}{\noindent{\bf Proof.\ }}
\newcommand{\qed}{\hfill $\square$ \bigskip}

\begin{document}
\title{Twin Vertices in Fault-Tolerant Metric Sets and Fault-Tolerant Metric Dimension of Multistage Interconnection Networks}

\author{
	S. Prabhu$^{1}$\thanks{Corresponding author: drsavariprabhu@gmail.com}
	\and
	V. Manimozhi$^{2}$
	\and
	M. Arulperumjothi$^{3}$
	\and
	Sandi Klav\v zar$^{4, 5, 6}$
}

\date{}

\maketitle
\begin{center}
$^{1}$  Department of Mathematics, Rajalakshmi Engineering College (Autonomous), Thandalam, Chennai 602105, India 
\medskip
	
$^{2}$ Department of Mathematics, Panimalar Engineering College, Chennai, India 600123	
\medskip

$^{3}$ Department of Mathematics, Saveetha Engineering College (Autonomous), Chennai 602105, India
\medskip

$^4$ Faculty of Mathematics and Physics, 
University of Ljubljana, Slovenia\\
\medskip

$^5$ Faculty of Natural Sciences and Mathematics,  University of Maribor, Slovenia\\
\medskip

$^6$ Institute of Mathematics, Physics and Mechanics, Ljubljana, Slovenia\\
\medskip

\end{center}

\begin{abstract}
\baselineskip16pt
A set of vertices $S\subseteq V(G)$ is a resolving set of a graph $G$ if for each $x,y\in V(G)$ there is a vertex $u\in S$ such that $d(x,u)\neq d(y,u)$. A resolving set $S$ is a fault-tolerant resolving set if $S\setminus \{x\}$ is a resolving set for every $x \in S$.  The fault-tolerant metric dimension (FTMD) $\beta'(G)$ of $G$ is the minimum cardinality of a fault-tolerant resolving set. It is shown that each twin vertex of $G$ belongs to every fault-tolerant resolving set of $G$. As a consequence, $\beta'(G) = n(G)$ iff each vertex of $G$ is a twin vertex, which  corrects a wrong characterization of graphs $G$ with $\beta'(G) = n(G)$ from [Mathematics 7(1) (2019) 78]. This FTMD problem is reinvestigated for Butterfly networks, Benes networks, and silicate networks. This extends partial results from [IEEE Access 8 (2020) 145435--145445], and at the same time, disproves related conjectures from the same paper. 
\end{abstract}



\baselineskip20pt
\newpage

\section{Introduction}

Let $G = (V(G), E(G))$ be a connected graph. The number of edges of a shortest $u,v$-path in $G$ is equal to the {\em distance} $d_G(u, v)$ (or simply $d(u, v)$ whenever $G$ is understood) between the vertices $u$ and $v$ of $G$. A finite collection of vertices $S\subseteq V(G)$ is a {\em resolving set} of $G$ if for every pair $a,b\in V(G)$ there is a $u\in S$ such that $d(a,u)\neq d(b,u)$. The \textit{metric dimension} $\beta(G)$ is the lowest cardinality of a resolving set of $G$.

This problem was initially investigated in the 1970s~\cite{HaMe76, Sl75}, and has been extensively investigated afterwards. The topic was first surveyed in 2003 by Chartrand and Zhang~\cite{ChZh03}. Very recently, a new excellent review on the metric dimension and its applications was written by Tillquist, Frongillo, and~Lladser~\cite{tillquists-2021+}. The survey includes $117$ references, we point here only to applications of resolving sets to verification of networks~\cite{BeEbEr06}, robot navigation~\cite{KhRaRo96}, and geometric routing protocols~ \cite{LiAb06}. 

Chartrand and Zhang in~\cite{ChZh03} also suggested using the members of a resolving set as sensors. Provided that, there will be a defective sensor leading to failure to acknowledge the burgler in the network. The idea arose from the assumption that a faulty sensor will not cause device failure because the remaining sensors will still be able to deal with the invasion force. It was formally introduced in~\cite{HeMoSl08}. A resolving set $F$ of a (connected) graph $G$ is a {\em fault-tolerant resolving set} if $F \smallsetminus \{u\}$ is a resolving set of $G$ for each $u \in F$. In other words, $F \subseteq G$ is a fault-tolerant resolving set if for every distinct $x,y \in V(G)$ there exist $u,v \in F$ such that $d(u,x) \neq d(u,y)$ and $d(v,x) \neq d(v,y)$. The least positive integer representing the cardinality of $F$ is $\beta'(G)$. The fault-tolerant metric dimension, like that of the metric dimension, has already been widely investigated, cf.~\cite{azeem-2021, basak-2020, JaSaCh09, raza-2019, SaJaCh14, SiBoMa18}.    

In~\cite{estrada-2015}, Estrada-Moreno, Rodr\'{\i}guez-Vel\'{a}zquez, and Yero introduced and studied the $k$-metric dimension, $k\ge 1$, as a common generalization of the metric dimension (the case $k=1$) and the fault-tolerant metric dimension (the case $k=2$).  The concept was later investigated in several papers, cf.~\cite{bailey-2019, corr-2021, estrada-2016}. We also emphasize that in~\cite{beardon-2019} the $k$-metric dimension was investigated on general metric spaces, while the fractional $k$-metric dimension is the core concept of~\cite{kang-2019}. In parallel to the survey~\cite{tillquists-2021+}, another excellent survey was posted by Kuziak and Yero~\cite{kuziak-2021}. The latter survey focuses on variants of the metric dimension, and in particular gives a thorough list on results on the $k$-metric dimension with a special attention to the case $k=2$, that is, to the fault-tolerant metric dimension. 

In the next section, we observe that each twin vertex of $G$ is a member of a fault-tolerant resolving set. As a consequence, we deduce that $\beta'(G) = n(G)$ iff each $v$ of $G$ is a twin. (Here and later, $n(G)$ denotes the order of $G$.) This result corrects a characterization with $\beta'(G) = n(G)$ that was incorrect in~\cite{RaHaIm19}. In Section~\ref{sec:butterfly-Benes}, we determine the FTMD of Butterfly, Benes, and silicate networks. For each of these networks, partial results were reported in~\cite{HaKhMa20}. Moreover, our results disprove related conjectures from the same paper. 

\section{Twin vertices}
\label{sec:twin}

The {\em open neighborhood} of a vertex $u\in V(G)$ is $N(u)=\{v\in V(G):uv\in E(G)\}$, the {\em closed neighborhood} of $u$ is $N[u]=N(u)\cup\{u\}$. Different vertices $u, v\in V(G)$ are {\em twins} if either $N[u]=N[v]$ or $N(u)=N(v)$. We further say that $u$ is a {\em twin vertex}, if there exists $v\ne u$ such that $u$ and $v$ are twins. Note that if $N[u]=N[v]$, then $uv\in E(G)$, and if $N(u)=N(v)$, then $uv\notin E(G)$. Twins are a very natural concept, hence no wonder that they were differently named in the literature. For instance, many authors say that $u$ and $v$ are {\em true twins} or {\em adjacent twins} when $N[u]=N[v]$ holds, and {\em false twins} or {\em non-adjacent twins} when $N(u)=N(v)$, cf.~\cite{barr-2019, lin-2018}. 

The role of twins for the metric dimension has been clarified in~\cite{HeMoPe10}. For the fault-tolerant metric dimension, we have the following observation. 

\begin{lemma}
\label{lem:all-twins}
If $u$ is a twin vertex of a graph $G$, and $S$ is a fault-tolerant resolving set of $G$,  then $u\in S$. In particular, if $S$ is the set of twin vertices of $G$, then $\beta'(G) \ge |S|$.
 \end{lemma}

\proof
Let $u$ be a twin vertex of $G$ and let $v\ne u$ be a vertex such that $u$ and $v$ are twins. Let $S$ is a fault-tolerant resolving set of $G$. Since $d(u,x) = d(v,x)$ for every $x\ne u,v$ (cf.~\cite[Lemma 2.3]{HeMoPe10}), and since $S$ is a resolving set, at least one of $u$ and  $v$ belongs to $S$. If $u\in S$ there is  nothing to prove. Assume hence that $v\in S$. Since $S$ is a fault-tolerant resolving set, the set $S\setminus \{v\}$ is a resolving set. If $u\notin S$, then the twins $u$ and $v$ have the same distance to all the vertices of the resolving set $S\setminus \{v\}$, a contradiction.  
\qed

Lemma~\ref{lem:all-twins} will be very useful in our consideration of interconnection networks in the next section. Before that, we point to a couple of errors from \cite{RaHaIm19}. 

Chartrand, Eroh, Johnson, and Oellermann~\cite{ChErJo00} proved that (i) $\beta(G) = n(G)-1$  iff $G$ is complete; and that (ii) if $n(G)\ge 4$, then $\beta(G) = n(G)-2 \Leftrightarrow$ $G$ belongs to one of the following graph classes $K_s + \overline{K_t}$ $(s \geq 1, t \geq 2)$, $K_ {s,t}$ $(s,t \geq 1)$, and $K_s + (K_1 \cup K_t)$ $(s,t \geq 1)$. Parallel results to (i) and (ii) were  erroneously  claimed in~\cite[Theorem 7]{RaHaIm19} and~\cite[Theorem 8]{RaHaIm19}, respectively. In particular, in the first of these results the authors claim that  $\beta'(G) = n(G)\Leftrightarrow G$ is complete. The correct result reads as follows. 
  
\begin{proposition} 
\label{prop:dim'=n}
Let $G$ be a connected graph. Then $\beta'(G) = n(G)$ if and only if each vertex of $G$ is a twin vertex.
\end{proposition}

\proof
If each vertex of $G$ is a twin vertex, then $\beta'(G) = n(G)$ by Lemma~\ref{lem:all-twins}. 

Conversely, suppose that $u\in V(G)$ is not a twin vertex. Then we claim that $V(G)\setminus \{u\}$ is a fault-tolerant resolving set. To show it, we only need to consider a pair of vertices $u$ and $x$, where $x$ is an arbitrary vertex from $V(G)\setminus \{u\}$. Clearly, $0 = d(x,x) < d(u,x)$. In addition, since $u$ and $x$ are not twins, there exits a vertex $y$ adjacent to exactly one of $u$ and $x$. Hence $d(u,y) \ne d(x,y)$. We conclude that if $G$ contains a vertex which is not a twin, then $\beta'(G) < n(G)$.
\qed

The graphs from Proposition~\ref{prop:dim'=n} can be described as follows. Its vertex set can be divided into disjoint parts of the order at least two, each part containing vertices that are pairwise twins. Moreover, each of these parts induces either a full graph or a null graph. Between each pair of these parts, there are either all possible edges or none at all. 

As mention in the introduction, the fault-tolerant metric dimension is just the $2$-metric dimension, and the two results of this section were actually already established in~\cite{estrada-2015}.  In particular, our Proposition~\ref{prop:dim'=n} is~\cite[Corollary 5]{estrada-2015}. We have nevertheless derived these results also here from the following reasons. First, in this way we are self-sufficient. Second, the arguments are fairly simple. And finally,~\cite[Corollary 5]{estrada-2015} is deduced from~\cite[Proposition 4]{estrada-2015}, where the latter result is stated in a slightly different and more involved language.  

\section{Butterfly, Benes, and silicate networks}
\label{sec:butterfly-Benes}

The illustration of a multistage interconnected model as graph has processors as vertices and interconnections between processors as edges. The structural features of interconnection networks were investigated in \cite{KrSn86}. Here, we determine the fault-tolerant metric dimension for two classes of such networks---butterfly networks and Benes networks. These networks were presented by Manuel et al.~\cite{MaAbRa08}, and investigated from differents angles, see~\cite{MaAbRa08} for their metric dimension, \cite{InHaMa14} for their degree-based topological indices, \cite{RaMaPa16} for their Wiener index, \cite{NuNaUd20} for their Zagreb indices and polynomials, and~\cite{MaRaRa13} for their crossing number.

\begin{figure}[ht!] 
	\centering
    \includegraphics[scale=.5]{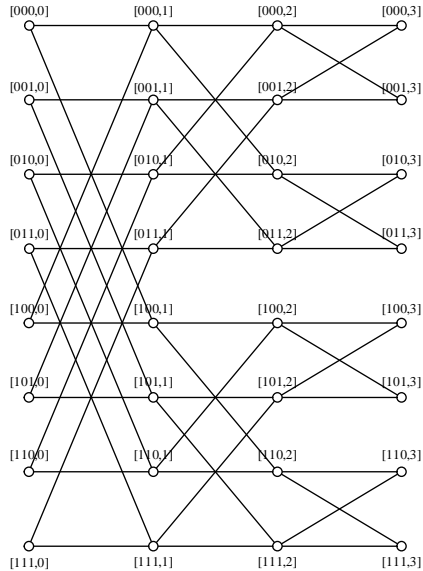}
	\caption {Normal representation of $BF(3)$} 
	\label{bfdefn}
\end{figure}

For $r\ge 3$, the $r$ dimensional {\em butterfly network} $BF(r)$ is defined as follows. Its vertices are pairs $[s, j]$, where $s$ runs over all $r$-bit binary strings, and $j\in \{0,1,\ldots, r\}$. The vertices $[s, j]$ and $[s ', j']$ are adjacent iff $|j-j'| = 1$, and either $s = s'$ or $s$ and $s'$ differ precisely in the $j^{\rm th}$ bit. Note that the order and the size of $BF(r)$ are $2^r (r + 1)$ and $r2^{r+1}$, respectively. In the {\em normal representation} of $BF(r)$, the first coordinate of the vertex is interpreted as the row of the vertex and its second coordinate is a column called {\em level} of the vertex. See Fig.~\ref{bfdefn} where $BF(3)$ is drawn. 

In~\cite[Corollary 5.3]{HaKhMa20} it was proved that $\beta'(BF(r))\leq 4\cdot 2^r$ and conjectured~\cite[Conjecture 5.4]{HaKhMa20}  that the equality holds here, that is, $\beta'(BF(r)) = 4\cdot 2^r$. In our next result we determine $\beta'(BF(r))$ which disproves the conjecture. 

\begin{theorem}
\label{thm:butterfly}
If $r \ge 3$, then $\beta'(BF(r))=2^{r+1}$.
\end{theorem}
  
\proof
From the definition of $BF(r)$ we infer that the vertices $[i-1 , 0]$ and $[i-1+2^{r-1}, 0]$  are twins for each $i\in [2^{r-1}]$. (We use the convention $[n] = \{1,\ldots, n\}$.) Moreover, for each $ i \in [2^{r-1}]$, the vertices $[2i-2 , r]$ and $[2i-1,r]$ are also twins. That is, each vertex from the first and the last level of $BF(r)$ is a twin vertex. Hence $BF(r)$ contains (at least) $2^{r+1}$ twin vertices, therefore $\beta'(BF(r)) \ge 2^{r+1}$ by Lemma~\ref{lem:all-twins}.

To prove that $\beta'(BF(r)) \le 2^{r+1}$, we are going to show that the set $X$ consisting of the vertices from the first and the last level of $BF(r)$ forms a fault-tolerant resolving set. For this sake, let $[s,i]$ and $[t,j]$ be arbitrary vertices of $BF(r)$. We need to show that they are distance distinguished by two vertices of $X$. If at least one of $[s,i]$ and $[t,j]$ lies in $X$, then the conclusion is clear, hence we may assume in the rest that $i, j\in [r-1]$. Suppose first that $i\ne j$ and assume w.l.g that $i < j$. Then $d([s,i], [s,0]) < d([t,j], [s,0])$ and $d([t,j], [t,r]) < d([s,i], [t,r])$. As $[s,0], [t,r]\in X$, we have required two vertices in this case. 

It remains to consider vertices $[s,i]$ and $[t,i]$, where $i\in [r-1]$ and $s\ne t$. Let $s=s_1\ldots s_r$ and $t=t_1\ldots t_r$. The following facts will be used in subsequent lines. If $x=x_1\ldots x_r$ and $d([x,0], [s,i]) = i$, then $x_{i+1} = s_{i+1}, \ldots,  x_{r} = s_{r}$. Similarly, if $d([x,r], [s,i]) = r-i$, then $x_{1} = s_{1}, \ldots,  x_{i} = s_{i}$. The binary labels $s$ and $t$ are different, hence they differ in at least one coordinate. Assume first that there exists an index $k > i$ such that $s_k \ne t_k$. Then we infer that $d([s,0],[s,i]) = i < d([s,0], [t,i])$ and that $d([t,0], [t,i]) = i < d([t,0], [s,i])$. Hence the vertices  $[s,0], [t,0]\in X$ distinguish $[s,i]$ and $[t,j]$. In the second case we may thus assume that $s_{i+1} = t_{i+1}, \ldots, s_{r} = t_{r}$. Let $k$ be an index such that  $s_k \ne t_k$. By the case assumption, $k\le i$. Consider now the vertices $[s,r]$ and $[t,r]$. Suppose that $d([s,r], [s,i]) = r-i = d([s,r], [t,j])$. Then, by the above fact, $t_{1} = s_{1}, \ldots,  t_{i} = s_{i}$. But this means that $t=s$, a contradiction. It follows that $d([s,r], [s,i]) = r-i < d([s,r], [t,i])$. Similarly we obtain that $d([t,r], [t,i]) = r-i < d([t,r], [s,i])$. Hence in this case the vertices  $[s,r]$ and $[t,r]$ from $X$ distinguish $[s,i]$ and $[t,j]$ and we are done. 
\qed

\begin{figure}[ht!] 
	\centering
     \includegraphics[scale=.5]{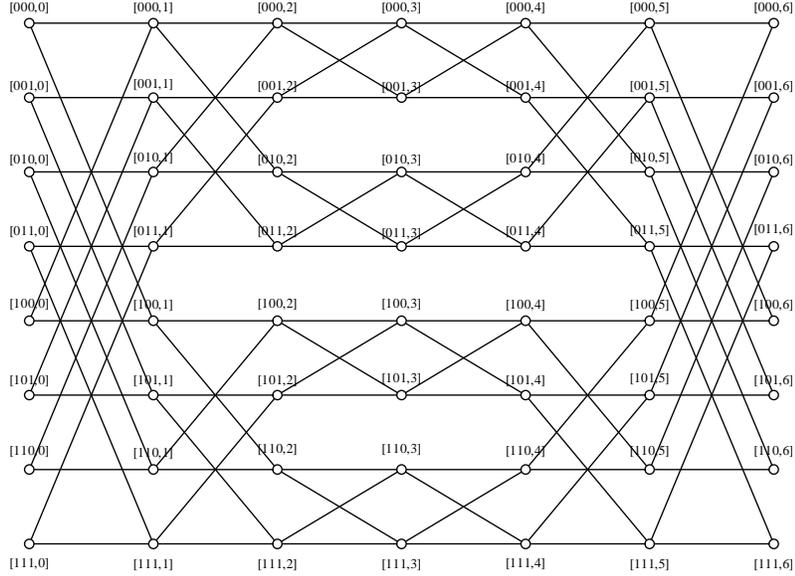} 
	\caption{ $B(3)$ drawn in its normal representation} 
	\label{bdefn}
\end{figure}

For $r\ge 3$, the $r$-dimensional {\em Benes network} $B(r)$ is defined similarly as $BF(r)$. The vertices are again ordered pairs $[s, i]$, where $s$ runs over all $r$-bit binary strings, but now $i\in \{0,1,\ldots, 2r\}$. 
The edges in $B(r)$ up to level $r$ are just as in $BF(r)$, while later on, the edges are vertically reflected. The formal definition should be clear from Fig.~\ref{bdefn}, where $B(3)$ is shown in its normal representation. 

The order and the size of $B(r)$ are $2^r (2r + 1)$ and $r2^{r+2}$, respectively.  In~\cite[Corollary 4.3]{HaKhMa20} it was proved that $\beta'(B(r))\leq 13\cdot 2^{r-1}$ and conjectured~\cite[Conjecture 4.4]{HaKhMa20}  that the equality holds here, that is, $\beta'(B(r)) = 13\cdot 2^{r-1}$. In our next result we determine $\beta'(B(r))$ which disproves the conjecture. 

\begin{theorem} \label{benes}
If $r\geq 3$, then $\beta'(B(r)) = 3\cdot2^r$.
 \end{theorem}

\proof
Following the arguments adopted in Theorem~\ref{thm:butterfly}, we observe first that each vertex from the first, the middle, and the last level of $B(r)$ is a twin vertex. Hence $B(r)$ contains (at least) $3\cdot 2^{r}$ twins, therefore $\beta'(B(r)) \ge 3\cdot 2^{r}$ by Lemma~\ref{lem:all-twins}.

To prove that $\beta'(B(r)) \le 3\cdot 2^{r}$, we claim that  the set $X$ consisting of the vertices from the first, the middle, and the last level of $B(r)$ forms a fault-tolerant resolving set. To show it we take arbitrary vertices $[s,i]$ and $[t,j]$ of $B(r)$. If at least one of them is in $X$, then two vertices from $X$ clearly distinguish them. Hence we may assume $i, j\in [2r-1]$. If $i\ne j$, where $i < j$, then $d([s,i], [s,0]) < d([t,j], [s,0])$ and $d([t,j], [t,2r]) < d([s,i], [t,2r])$. Hence it remains to consider vertices $[s,i]$ and $[t,i]$, where $i\in [2r-1]$ and $s\ne t$. If $i\in [r-1]$, then we proceed analogously as in the proof of Theorem~\ref{thm:butterfly}, hence we omit the details. And if $i\in  \{r+1, \ldots, 2r-1\}$, then by the vertical symmetry of $B(r)$ over the level $r$ (with respect to the normal representation of $B(r)$) we get the required conclusion is an analogous way. 
\qed

For Benes and butterfly networks, it was also conjectured in~\cite[Conjecture 7.1]{HaKhMa20} that the FTMD problem is polynomially solvable. Theorems~\ref{thm:butterfly} and~\ref{benes} clearly confirm the conjecture. 

\begin{figure}[ht!] 
 	\centering
 	\includegraphics[scale=.5]{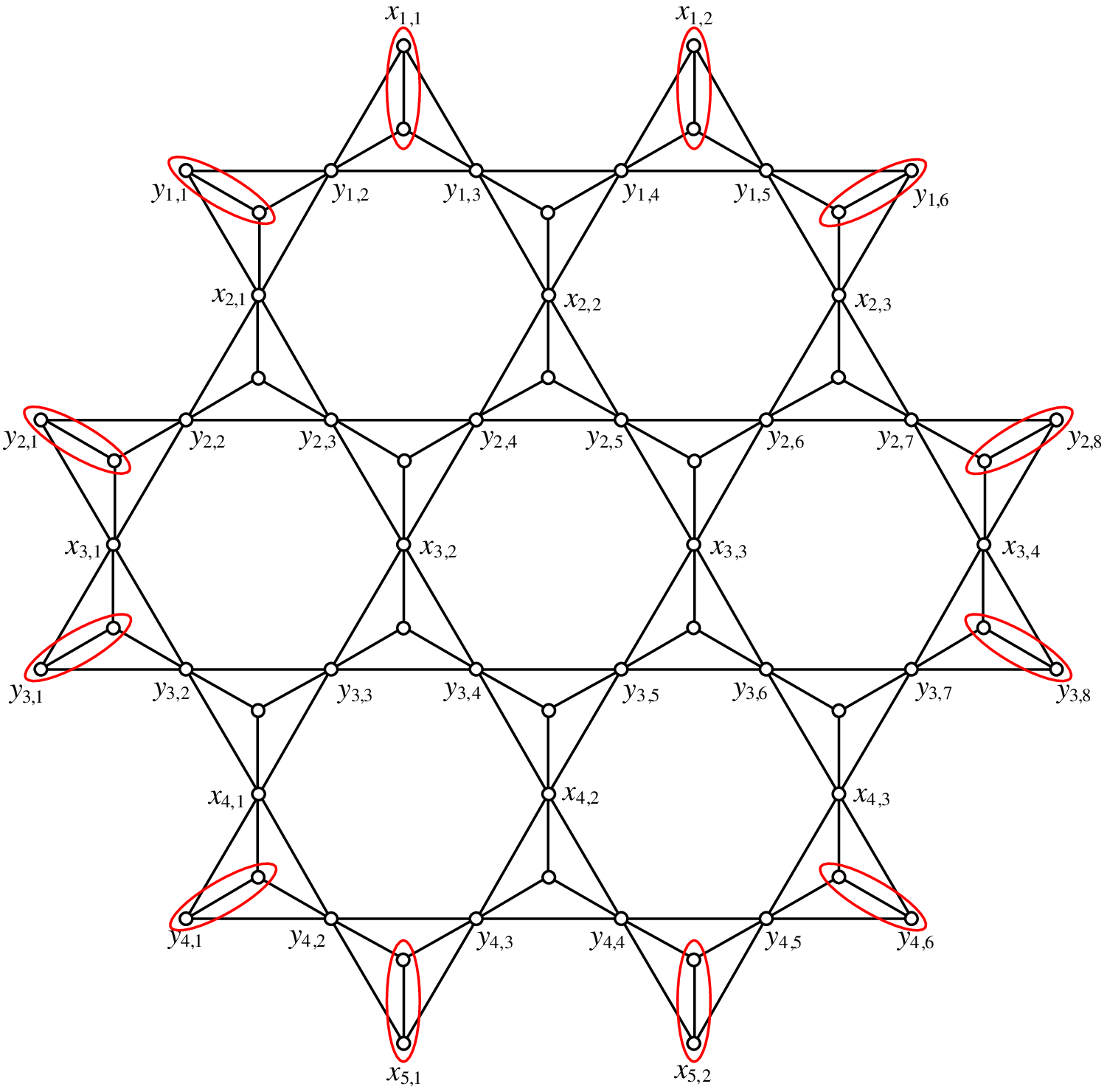} 
 	\caption {$SL(2)$ and its twins} \label{fig:slts}
 \end{figure}

The last networks we consider are the silicate networks $SL(n)$, $n>1$, which can be defined as follows. Begin with the $n^{\rm th}$ member of the circumcoronene homologous series $H_n$, which consists of the central hexagon and $n-1$ layers of hexagons around it. (Metric-based resolvability of these polycyclic aromatic hydrocarbons were very recently investigated in detail in~\cite{azeem-2021}.) Then subdivide each edge of it, and finally append a $K_4$ around the original vertices of $H_n$, where for each outer vertex of $H_n$ we need to add an additional vertex to accomplish $K_4$. The construction should be clear from Fig.~\ref{fig:slts}, where $SL(2)$ is drawn. 

More generally, the anionic substructure of silicates comprises (SiO$_4$) tetrahedra, which can crosslink by sharing common corners. Their primary variety has been widely explored, starting from mineralogy and topography to chemical, physical, and computer sciences. From a synthetic perspective, the tetrahedron's focal hub represents the silicon particle, while the corner ions speak to the oxygen atoms, while the corner vertices represent oxygen atoms. We get distinct silicate structures by arranging different ways of tetrahedra. Papers~\cite{HaIm14, LiWaWa17, MaRa09} derive some topological features of the silicate networks. This network family is a multi-stage interconnection network comparable to the hexagonal \cite{ChShKa90} and honeycomb \cite{St97} networks used in computer science networks~\cite{St97}. This is achieved by suggesting an addressing scheme to the vertices of the silicate networks as suggested for hexagonal in~\cite{GaNoSt02} and honeycomb~\cite{St97}. 

The metric dimension of the silicate networks $SL(n)$ was determined in 2011 as follows.

\begin{theorem} {\rm \cite{MaRa11}} 
If $n\ge 2$, then $\beta(SL(n)) = 6n$.
\end{theorem}

For the FTMD of silicate networks it was proved in~\cite[Proposition 6.4]{HaKhMa20} that  
$$6n+1 \leq \beta'(SL(n)) \leq 21n\,.$$
In our last result, we give an exact formula for FTMD of silicate networks.

\begin{theorem}
If $n\ge 2$, then $\beta'(SL(n))=12n$. 
 \end{theorem}
 
\proof
$SL(n)$ contains $6n$ pairs of twin vertices, as shown in Fig~\ref{fig:slts} for $SL(2)$. Hence $SL(n)$ contains $12n$ twin vertices and thus $\beta'(SL(n)) \ge 12n$ by Lemma~\ref{lem:all-twins}. To prove the reverse inequality, we claim that these $12n$ vertices form a fault-tolerant resolving set. Instead of proving it directly, we recall from the proof of~\cite[Theorem 5]{MaRa11} that the set of boundary vertices (the $6n$ outermost vertices in the standard drawing of $SL(n)$, cf.~Fig.~\ref{fig:slts}) of $SL(n)$ is a (minimum) resolving set of $SL(n)$. But then it straightforwardly follows that these boundary $6n$ vertices, together with their respective twins, form a (minimum) fault-tolerant resolving set.  
\qed

\section*{Acknowledgements}

Sandi Klav\v{z}ar acknowledges the financial support from the Slovenian Research Agency (research core funding P1-0297, and projects N1-0095, J1-1693, J1-2452).

\end{document}